\documentclass[11pt]{article}
\usepackage{latexsym}
\usepackage[all]{xy}
\usepackage{amsmath}
\usepackage{amssymb}
\usepackage{amsfonts}
\usepackage[applemac]{inputenc}
\RequirePackage{hyperref}

\newtheorem{thm}{Theorem}[section]
\newtheorem{prop}[thm]{Proposition}

\newtheorem{df}[thm]{Definition}
\newtheorem{cor}[thm]{Corollary}
\newtheorem{conj}{Conjecture}
\newtheorem{rmk}[thm]{Remark}

\title{\textbf{Some elementary remarks on lci algebraic cycles}}

\date{Notes for students\\
August 2016}

\author{M. Maggesi and G. Vezzosi \\
\small{Dipartimento di Matematica ed Informatica}\\ \small{Universit\`a di Firenze - Italy}}

\begin{document}

\maketitle

\begin{abstract}
\noindent In this short note, we simply collect some known results about representing algebraic cycles by various kind of ``nice'' (e.g.~smooth, local complete intersection, products of local complete intersection)  algebraic cycles, up to rational equivalence. We also add a few elementary and easy observations on these representation problems that we were not able to  locate in the literature. 
\end{abstract}


\section*{Introduction}

The question whether any algebraic cycle on a complex smooth (quasi-)projective variety $X$ is rationally equivalent to (a multiple of) an algebraic cycle belonging to some special and nice class, is an old problem. To our knowledge, it was first addressed by Hironaka \cite{hi} and Kleiman \cite{kl} who considered the class of \emph{smooth} cycles, by analogy with similar problems in topology originating form the Steenrod problem and from sebsequent work and questions by Borel-Haefliger. Their results, together with a subsequent relevant result by Hartshorne-Rees-Thomas \cite{hrt}, are briefly reviewed in Section \ref{sm} below.
One might ask the same questions for larger classes of algebraic cycles models. A natural choice, from the standpoint of derived algebraic geometry, is the class of \emph{local complete intersection} (lci) cycles. We collect below a few elementary results for this class of models, and more generally we discuss and compare three different classes of models: smooth, lci, and products of lci cycles. Though we were not able to make any sensitive progress in either of these cases, we decided to write this note because we think these representation problems for algebraic cycles are extremely interesting, though a bit neglected in the current research panorama on algebraic cycles.\\

Most, if not all, of the facts below hold for a complex quasi-projective smooth scheme (and often, over an arbitrary algebraically closed field).
However, we will stick to smooth complex projective varieties because these are our main objects of interest in this note.\\

\noindent  \textbf{Acknowledgments}. After obtaining our remarks (Section \ref{lcisection}), we discovered the following nice MathOverflow page
http://mathoverflow.net/questions/60434 where similar questions are discussed. Though few details are given, we think that Sasha's idea on the proof of Proposition \ref{prelim} must have been very similar to ours, and he should be credited for this. The rest of the arguments in the above MO page, especially the one claiming that every Chern class (not only top ones) is in fact an lci cycle, are more obscure to us. As Burt Totaro pointed out to us, determinantal varieties are typically not lci, so the corresponding MO page argument is probably incomplete.\\
GV wishes to thank Burt Totaro for several useful e-mail exchanges related to the subject matter of this note, that revived his (i.e.~GV's) interest in representation problems for algebraic cycles, and corrected some mistakes. GV also wishes to thank Mauro Porta for many inspiring conversations, and some joint elementary experiments about the existence of a derived version of the Chow group.\\

\noindent  \textbf{Notations}. $X$ will be a complex smooth projective variety, and we will write $\mathrm{CH}_*(X)$ (respectively $\mathrm{CH}^*(X)$) for the Chow group (respectively, the Chow ring) of algebraic cycles modulo rational equivalence on $X$, \emph{with $\mathbb{Q}$-coefficients}. The word \emph{variety} is used as in Fulton's \cite{fu} to mean a reduced, irreducible scheme of finite type over the base field (which will be $\mathbb{C}$ in this note). For $Z\subset X$ a closed subscheme, $[Z]$ will denote the algebraic cycle (modulo rational equivalence) associated to $Z$, as in \cite[1.5]{fu}.\\

\section{Smooth, lci, and plci algebraic cycles}

\begin{df} Let $X$ be a smooth complex projective variety.
\begin{itemize}
\item An algebraic cycle $\mathfrak{z}$ on $X$ is said to be \emph{lci} or \emph{local complete intersection} (respectively, \emph{smooth}) if there exist closed subschemes $i_j: Z_j \hookrightarrow X$,  and $\alpha_j \in \mathbb{Q}$, $j=1, \ldots, n$, such that each $i_j$ (respectively, each $Z_j$) is a regular closed immersion\footnote{Or, equivalently, an lci morphism.} (resp. is smooth over $\mathbb{C}$), and $\mathfrak{z}$ is rationally equivalent to the algebraic cycle $\sum^n_{j=1}\alpha_j [Z_j]$ (i.e.~$[\mathfrak{z}]= \sum^n_{j=1}\alpha_j [Z_j]$ in $\mathrm{CH}_{*}(X)$). We denote by $\mathrm{CH}^{\mathrm{lci}}_{*}(X)$ (respectively, $\mathrm{CH}^{\mathrm{sm}}_{*}(X)$)  the sub-$\mathbb{Q}$-vector space of $\mathrm{CH}_{*}(X)$ generated by lci (resp., smooth\footnote{This definition of smooth algebraic cycles is equivalent to the one considered in \cite{hi, kl, hrt}, since a smooth $\mathbb{C}$-scheme has a finite number of irreducible components, each one being again smooth over $\mathbb{C}$.}) cycles.
\item An algebraic cycle $\mathfrak{z}$ on $X$ is said to be \emph{plci} (or \emph{product of local complete intersection}) if it is rationally equivalent to a product of lci cycles inside the graded $\mathbb{Q}$-algebra $\mathrm{CH}^*(X)$. We denote by $\mathrm{CH}^{*}_{\mathrm{plci}}(X)$ the sub-$\mathbb{Q}$-vector space of $\mathrm{CH}^{*}(X)$ generated by lci cycles (or, equivalently, the sub-$\mathbb{Q}$-algebra of $\mathrm{CH}^{*}(X)$ generated by lci cycles).
\end{itemize} 
\end{df}

It is natural to formulate the following questions (that we state as conjectures because logical implications between them will be easier to state).

\begin{conj}\label{smooth}
$\,\,\,\,\,\, \mathrm{CH}^{\mathrm{sm}}_{*}(X) = \mathrm{CH}_{*}(X).$
\end{conj}

\begin{conj}\label{lci}
$\,\,\,\,\,\, \mathrm{CH}^{\mathrm{lci}}_{*}(X) = \mathrm{CH}_{*}(X).$
\end{conj}

\begin{conj}\label{plci}
$\,\,\,\,\,\, \mathrm{CH}_{\mathrm{plci}}^{*}(X) = \mathrm{CH}_{*}(X).$
\end{conj}

We have obvious implications $$\mathrm{Conjecture\,\, \ref{smooth}} \Rightarrow \mathrm{Conjecture\,\, \ref{lci}} \Rightarrow \mathrm{Conjecture \,\, \ref{plci}}.$$

To our knowledge, none of the previous implications are known to be reversable in general.\\ 

\begin{rmk} \emph{By Bertini's theorem, \emph{global} complete intersections (gci) in $X= \mathbb{P}^n$ or in  a smooth affine variety $X= \mathrm{Spec}\, A$ are smoothable, i.e. can be deformed to smooth subschemes. Moreover, the parameter space of the deformation can always be assumed to be a non-singular curve (see e.g. the argument right before \cite[Example 29.0.1]{hadef}).
For such $X$, we have $\mathrm{CH}^{\mathrm{gci}}_{\mathrm{alg},\, *}(X)= \mathrm{CH}^{\mathrm{sm}}_{\mathrm{alg},\,*}(X)$, for the Chow subgroups of gci and smooth cycles, modulo \emph{algebraic equivalence}. However, on a general smooth projective $X$, it is not true that an lci subscheme is smoothable: a singular $(-1)$-curve on a smooth surface is rigid, hence is not smoothable. Further examples of non-smoothable lci curves in $\mathbb{P}^3$ are given in \cite[\S \, 5]{hahi}. Since we are interested in algebraic cycles, and not in subschemes, it might nonetheless be true that $\mathrm{CH}^{\mathrm{lci}}_{\mathrm{alg},\, *}(X)= \mathrm{CH}^{\mathrm{sm}}_{\mathrm{alg},\,*}(X)$, or that $\mathrm{CH}^{\mathrm{lci}}_{ *}(X)= \mathrm{CH}^{\mathrm{sm}}_{*}(X)$\footnote{Another reason to believe this comes from the results of Hironaka and Kleiman, see Section \ref{sm}}. For example, it is possible that an arbitrary given lci subscheme becomes smoothable after adding a suitable smooth subscheme (this is indeed known to be the case in codimension $1$, e.g. for the $(-1)$ curves mentioned above).
This remark is essentially due to Burt Totaro, who also provided the reference \cite{hahi}. We wish to thank him for his clarifications.}
\end{rmk}

\subsection{Discussion of Conjecture \ref{smooth}} \label{sm}

Conjecture \ref{smooth} have been addressed by Hironaka (for Chow groups with integer and rational coefficients, \cite{hi}), and by Kleiman (for Chow groups with rational coefficients, \cite{kl}). Kleiman proved that for a connected, smooth (quasi-)projective scheme $Y$ of dimension $d$ over an algebraically closed field, then $\mathrm{CH}_i^{\mathrm{sm}}(Y) = \mathrm{CH}_i(Y)$, for any $i < (d+2)/2$ (\cite[Theorem 5.8]{kl}). Later,  Hartshorne-Rees-Thomas proved that Conjecture \ref{smooth} is false \emph{with integer coefficients}, by producing algebraic cycles with integer coefficients on $\mathrm{Grass}_{3,3}$ that are not smoothable even up to homological equivalence (\cite[Theorem 1]{hrt}).  However, the examples in \cite{hrt} are not known to contradict Conjecture \ref{smooth}, i.e.\ it is not known whether some integer multiples of them are smoothable.\\

The next remark shows that, in order to check whether Conjecture \ref{smooth} is true for $X$, it is enough to test it on all Chern classes of all vector bundles on $X$.
 
\begin{prop} Let $X$ be a smooth complex projective variety. If for all vector bundles $E$ on $X$, we know that $c_i(E)$ is a smooth cycle, for any $i\geq0$, then Conjecture \ref{smooth} is true.
\end{prop}

\noindent \textbf{Proof.} Using the syzygy theorem as in the proof of \cite[Lemma 5.4]{kl}, we know that, for any closed irreducible codimension $p$ subvariety $Z \subset X$, there is a vector bundle $E$ on $X$, and an integer $n$, such that the following equality $$\pm (p-1)! [Z]= c_{p}(E) - n\mathrm{H}^p$$ holds in $\mathrm{CH}^*(X)$. Here $\mathrm{H}$ is the class of an hyperplane section, hence $\mathrm{H}^p$ is smooth for any $p$, and we conclude.
\hfill $\Box$ \\

\subsection{Discussion of Conjecture \ref{lci}}\label{lcisection}

First of all, as in the case of smooth cycles, in order to prove Conjecture \ref{lci} for $X$, it is enough to test it on all Chern classes of vector bundles on $X$:
 
\begin{prop} Let $X$ be a smooth complex projective variety. If for all vector bundles $E$ on $X$, we know that $c_i(E)$ is an lci cycle, for any $i\geq0$, then Conjecture \ref{lci} is true.
\end{prop}

\noindent \textbf{Proof.} Using the syzygy theorem as in the proof of \cite[Lemma 5.4]{kl}, we know that, for any closed irreducible codimension $p$ subvariety $Z \subset X$, there is a vector bundle $E$ on $X$, and an integer $n$, such that the following equality $$\pm (p-1)! [Z]= c_{p}(E) - n\mathrm{H}^p$$ holds in $\mathrm{CH}^*(X)$. Here $\mathrm{H}=c_1(\mathcal{O}(1))$ is the class of an hyperplane section, hence $\mathrm{H}^p$ is lci for any $p$ (e.g.\ by Corollary \ref{ctop} below), and we conclude.
\hfill $\Box$ \\

Therefore we can reduce to study the property of being lci for Chern classes of vector bundles on $X$. For the first Chern classes, it is classical that they all belong to $\mathrm{CH}_{d-1}^{\mathrm{sm}}(X)$, hence to $\mathrm{CH}_{d-1}^{\mathrm{lci}}(X)$, for $d= \dim \, X$. A further observation in this direction will be a corollary of the following

\begin{prop}\label{prelim}
  \label{prop:vandermonde}
  Let $E$ be a rank $ r$ vector bundle on $X$, and $L$ be a globally generated line bundle on $X$, such that $E \otimes L$ is globally generated. Then, the Chern classes \( c_1(L)^ic_{r-i}(E \otimes L), \)
  are lci for all $i = 0,\dots, r$.
\end{prop}
\noindent \textbf{Proof.}
  Define $x_i := c_1(L)^ic_{r-i}(E\otimes L)$.  By developing the Chern
  class of the tensor product
  $E\otimes L^{\otimes \, m+1} = (E\otimes L) \otimes L^{\otimes \, m}$ we get
  \begin{equation*}
    c_r((E\otimes L)\otimes L^{\otimes \, m}) 
    = \sum_{i=0}^r c_1(L^{\otimes \, m})^i c_{r-i}(E\otimes L) 
    = \sum_{i=0}^r m^i x_i
  \end{equation*}
  where we have used that $c_{1}(L^{\otimes \, m})= m c_1 (L) $.
  By taking $m = 0,\dots, r$ we get a system of $(r+1)$ linear equations
  in the $(r+1)$ unknowns $x_0, \dots, x_r$ with integer coefficients, and
  constant terms given by the \ classes $c_r(E\otimes L^{m+1})$, for $m= 0,\dots, r$.
  Note that the top Chern class $c_{top}(F)$ of any globally generated vector bundle $F$ on $X$
  is represented by the zero scheme $Z(s)$ of some global section of $F$, and $Z(s)\hookrightarrow X$ is a regular closed immersion of codimension $\mathrm{rk} (F)$, since $X$ is Cohen-Macaulay (see e.g.~\cite[Thm. 8.21A (c)]{habook}). Hence $c_{top}(F)$ is an lci cycle for any globally generated vector bundle $F$ on $X$; hence this is true for $c_r(E\otimes L^{m+1})$, for $m= 0,\dots, r$, by our hypothesis on $L$.
  
  The matrix of the corresponding homogeneous system is of the form
  \begin{equation*}
    \begin{bmatrix}
      1 & 0 & 0 &\dots & 0 \\
      1 & 1 & 1 & \dots & 1 \\
      1 & 2 & 4 & \dots & 2^r \\
      \vdots & & \vdots & & \vdots \\
      1 & r & r^2 &\dots & r^r \\
    \end{bmatrix}
  \end{equation*}
  which has a non-vanishing determinant since it is a Vandermonde
  matrix whose rows are all different.  By solving the system, we find
  $x_0,\dots,x_r$ expressed as rational linear combination of the
  lci\ classes $c_r(E\otimes L),\dots,c_r(E\otimes L^{\otimes \, r+1})$.  Thus
  all the $x_i$'s are lci.\footnote{By observing that
    $x_0=c_r(E\otimes L)$ and $x_r=c_1(L)^r$ are lci, we could
    reduce ourself to consider a smaller linear system of $r-1$
    equations in $r-1$ unknowns.  However, this does not yield any
    significant simplification in our argument.}
\hfill $\Box$ \\

\begin{cor}\label{ctop} Let $E$ be a vector bundle of rank $r$ on a smooth projective variety $X$. 
  The top Chern class $c_r(E)$ is lci.
\end{cor}
\noindent \textbf{Proof.} Pick a line bundle $L$ on $X$ as in the statement of Proposition \ref{prelim} (it is classical that such a sufficiently ample line bundle exists).
  Consider the top Chern class of $E\otimes L$;  we have
  \begin{equation*}
    c_r(E\otimes L) = \sum_{i=0}^r c_1(L)^ic_{r-i}(E).
  \end{equation*}
  Thus
  \begin{equation*}
    c_r(E) = c_r(E\otimes L) - \sum_{i=1}^rc_1(L)c_{r-i}(E).
  \end{equation*}
  By Proposition \ref{prelim}, the right-hand side of
  the previous equation is indeed a sum of lci classes.
\hfill $\Box$ \\

At the moment we are unable to prove that $c_i(E)$ is lci for any vector bundle $E$ on $X$, and $1< i < \mathrm{rank}\, E$ (i.e we are not able to prove Conjecture \ref{lci}).


\begin{rmk} \emph{Let us compare Corollary \ref{ctop} with the best, up to now, smoothing result for cycles with rational coefficients, i.e.~\cite[Theorem 5.8]{kl}. Let $d:= \dim \, X$. As already recalled above, 
Kleiman shows that $\mathrm{CH}_{i}(X)=\mathrm{CH}^{\mathrm{sm}}_{i}(X)$ for $i < (d +2)/2$. This result is much stronger (and effective) than Corollary \ref{ctop}, for $\mathrm{rank}\, E > (d-2)/2$. However, for $2 \leq\mathrm{rk}\, E \leq (d-2)/2$, \cite[Theorem 5.8]{kl} does not cover Corollary \ref{ctop}. }
\end{rmk}


\subsection{Discussion of Conjecture \ref{plci}} 

Conjecture \ref{plci} is the weakest one among the three we have listed. We remark the following
\begin{prop} If $c_i(E)$ is plci, for any vector bundle $E$ over $X$, and for any $i \geq 0$, then Conjecture \ref{plci} is true.
\end{prop}

\noindent \textbf{Proof.} 
Simply recall that the Chern character $\mathrm{Ch}: K_0(X) \otimes \mathbb{Q} \to \mathrm{CH}^*(X)$ is an isomorphism of $\mathbb{Q}$-algebras, since $X$ is smooth, and that the Chern character is a polynomial with rational coefficients in the Chern classes. 
\hfill $\Box$ \\

Suppose we know that Conjecture \ref{plci} is true. Then it would be tempting to deduce Conjecture \ref{lci} by a deformation argument saying that the intesection of two lci cycles is rationally equivalent to an lci cycle. In order to do this, one might use (at least) two strategies: deformation to the normal cone or a moving lemma. We have not fully explored the first approach, but we are a bit skeptical about the second one because the classical moving lemma deforms cycles, not specified subvarieties (or rather, it deforms a specified subvariety into a full cycle). So it is not obvious to us how a moving lemma argument might be combined with a deformation theoretic argument on the moduli space of lci subschemes in $X$, in order to prove that the product of any two lci cycles is indeed an lci cycle (i.e. that plci $\Rightarrow$ lci).

\end{document}